\documentstyle[12pt]{article}
\makeatletter
\renewcommand\thesection{\@Roman\c@section}
\renewcommand\thesubsection{\thesection.\@arabic\c@subsection}
\makeatother

\begin{document}
\begin{titlepage}
\begin{flushright}
math.QA/9812084
\end{flushright}
\vskip.3in

\begin{center}
{\Large \bf Level-One Representations and Vertex Operators of Quantum Affine 
Superalgebra $U_q[\widehat{gl(N|N)}]$}
\vskip.3in
{\large Yao-Zhong Zhang}
\vskip.2in
{\em Department of Mathematics, University of Queensland, Brisbane,
     Qld 4072, Australia

Email: yzz@maths.uq.edu.au}
\end{center}

\vskip 2cm
\begin{center}
{\bf Abstract}
\end{center}

Level-one representations of the quantum affine superalgebra
$U_q[\widehat{gl(N|N)}]$ associated to the appropriate
non-standard system of simple roots and $q$-vertex operators 
(intertwining operators) associated with the level-one modules
are constructed explicitly in terms of free bosonic fields.


\end{titlepage}


\def\a{\alpha}
\def\b{\beta}
\def\d{\delta}
\def\e{\epsilon}
\def\ve{\varepsilon}
\def\g{\gamma}
\def\k{\kappa}
\def\l{\lambda}
\def\o{\omega}
\def\t{\theta}
\def\s{\sigma}
\def\D{\Delta}
\def\L{\Lambda}

\def\G{{gl(N|N)}}
\def\hG{{\widehat{gl(N|N)}}}
\def\R{{\cal R}}
\def\hR{{\hat{\cal R}}}
\def\C{{\bf C}}
\def\P{{\bf P}}
\def\Z2{{{\bf Z}_2}}
\def\T{{\cal T}}
\def\H{{\cal H}}
\def\trho{{\tilde{\rho}}}
\def\tphi{{\tilde{\phi}}}
\def\tT{{\tilde{\cal T}}}
\def\uqgnh{{U_q[\widehat{gl(N|N)}]}}
\def\uqg1h{{U_q[\widehat{gl(1|1)}]}}


\def\beq{\begin{equation}}
\def\eeq{\end{equation}}
\def\bea{\begin{eqnarray}}
\def\eea{\end{eqnarray}}
\def\ba{\begin{array}}
\def\ea{\end{array}}
\def\no{\nonumber}
\def\lt{\left}
\def\rt{\right}
\newcommand{\bq}{\begin{quote}}
\newcommand{\eq}{\end{quote}}

\newtheorem{Theorem}{Theorem}
\newtheorem{Definition}{Definition}
\newtheorem{Proposition}{Proposition}
\newtheorem{Lemma}{Lemma}
\newtheorem{Corollary}{Corollary}
\newcommand{\proof}[1]{{\bf Proof. }
        #1\begin{flushright}$\Box$\end{flushright}}

\newcommand{\sect}[1]{\setcounter{equation}{0}\section{#1}}
\renewcommand{\theequation}{\thesection.\arabic{equation}}

\sect{Introduction\label{intro}}

The algebraic analysis approach \cite{Dav93,Jim94}
based on quantum affine algebra symmetries
enables one not only to solve 
massive or off-critical integrable models directly in the
thermodynamic limit but also to compute their
correlation functions \cite{Kor93} and form factors \cite{Smi92}
in the form of integrals by applying the techniques similar to those
used so successfully in the critical cases (see, c.f. \cite{Itz88}). 
The key components behind this method are
infinite dimensional highest weight representations of the
quantum affine algebras and the
corresponding $q$-vertex operators \cite{Fre92} which are intertwiners of these
representations. As in the critical cases, this procedure requires
the explicit construction of the highest weight
representations and vertex operators in terms of free bosonic fields.

By now, the level-one representations and vertex operators have been
constructed in terms of free bosons for most quantum
affine bosonic algebras (see, e.g. \cite{Fre88,Ber89,Awa94,Koy94,Jin95,Jin97}).
In contrast, much less has been known for the case of quantum affine
superalgebras. For the type I quantum affine superalgebra 
$U_q[\widehat{gl(M|N)}]$, $M\neq N$, the level-one representations
and vertex operators have been investigated in \cite{Kim97} (see
\cite{Awa97} for a level-$k$ free boson realization of
$U_q[\widehat{sl(2|1)}]$). In particular, the level-one irreducible
highest weight representations of $U_q[\widehat{gl(2|1)}]$ were
studied in some details and the correpsonding characters were derived
\cite{Kim97}. These representations have been re-examined and used 
to compute the correlation functions of the $q$-deformed
supersymmetric $t$-$J$ model in \cite{Yan98}.

So far in the literature, the very interesting case of $M=N$ has 
been largely ignored. The only exception is
\cite{Gad97} where the special case of $M=N=2$ was treated and 
the type I vertex operators involving infinite
dimensional evaluation (or level-zero) representations were also constructed
for this special case. By contrast, we shall consider the general
$M=N$ case and investigate both type I and type II 
vertex operators with respect to finite
dimensional evaluation modules.
The $M=N$ case is interesting since it seems to us that
$U_q[\widehat{gl(N|N)}]$ is the only untwisted superalgebra which has
a non-standard system where all simple roots are odd or fermionic.
It also seems to be the only superalgebra where a vertex type quasi-Hopf
twistor can be constructed \cite{Zha98} and thus the correpsonding
elliptic quantum supergroup ${\cal A}_{q,p}[\widehat{gl(N|N)}]$
can be introduced. 

In this paper, we construct a level-one representation of $\uqgnh$
by bosonizing the Drinfeld generators. We also construct the
vertex operators associated with the level-one representations in terms
of the free bosonic fields. 

The layout of this paper is the following. In section 2, we describe
the Drinfeld realization \cite{Dri88} of $\uqgnh$ in the non-standard
system of simple roots and determine the ``main terms" \cite{Cha91} in
the coproduct formulae of the Drinfeld generators. In section 3, we derive
the $2N$-dimensional evaluation (or level-zero) representations of
$\uqgnh$. In section 4, we investigate the bosonization of $\uqgnh$ and
construct an explicit level-one representation in terms of free bosonic
fields. Section 5 is devoted to the study of the bosonization of the
level-one vertex operators.

\sect{Quantum Affine Superalgebra $\uqgnh$}

As is well-known, a given Kac-Moody superalgebra
\cite{Kac78} allows many inequivalent systems of simple roots.
A system of simple roots is called distinguished if it has
minimal odd roots. Let $\{\a_i,~i=0,1,\cdots,2N-1\}$ denote a
chosen set of simple roots of the affine superalgebra
$\hG$. Let $(~,~)$ be a fixed 
invariant bilinear form on the root space. Let $\H$ be
the Cartan subalgebra and throughout we identify the
dual $\H^*$ with $\H$ via $(~,~)$. As is shown in \cite{Zha98},
$\hG$ has a simple root system in which
all simple roots are odd (or fermionic). This system can be constructed
from the distinguished
simple root system by using the ``extended"
Weyl operation \cite{Fra89} repeatedly. We have the following simple
roots, all of which are odd (or fermionic)
\bea
&&\a_0=\d-\ve_1+\ve_{2N}\,,\no\\
&&\a_l=\ve_l-\ve_{l+1}\,,~~~~~l=1, 2,\cdots,2N-1\label{roots}
\eea
with $\d,~\{\ve_k\}_{k=1}^{2N}$ satisfying
\beq
(\d,\d)=(\d,\ve_k)=0,~~~~(\ve_k,\ve_{k'})=(-1)^{k+1}\d_{kk'}.
\eeq
Such a simple root system is usually called non-standard. 
The generalized symmetric Cartan matrix $(a_{ii'})$ takes the form
\bea
&&a_{01}\equiv (\a_0,\a_1)=-1,~~~~a_{0,2N-1}\equiv (\a_0,\a_{2N-1})
  =1,\no\\
&&a_{ll'}\equiv (\a_l,\a_{l'})=(-1)^{l+1}(\d_{l,l'-1}-\d_{l,l'+1}),~~~
  l,l'=1,2,\cdots,2N-1.
\eea
This Cartan matrix is degenrate. To obtain a non-degenerate Cartan
matrix, we extend \cite{Kho91} $\H$
by adding to it the element
\beq
\a_{2N}=\sum_{k=1}^{2N}\ve_k\label{root-2n}.
\eeq
In the following, we denote by $\tilde{\H}$ the extended Cartan
subalgebra and by $\tilde{\H}^*$ the dual of $\tilde{\H}$.
The enlarged Cartan matrix has the following extra matrix elements:
\beq
a_{2N,2N}\equiv (\a_{2N},\a_{2N})=0,~~~
a_{i,2N}\equiv (\a_i,\a_{2N})=2 \cdot(-1)^{i+1}.
\eeq

Let $\{h_0,h_1,\cdots,h_{2N},d\}$ be a basis of $\tilde{\cal H}$, where 
$h_{2N}$ is the element in $\tilde{\H}$ corresponding to $\a_{2N}$ and $d$
is the usual derivation operator.
We shall write $h_i=\a_i~(i=0,1,\cdots,2N)$ with $\a_i$ given by
(\ref{roots}, \ref{root-2n}). Let $\{\L_0,\L_1,\cdots,\L_{2N},c\}$
be the dual basis with $\L_j$ being fundamental weights and $c$ the
canonical central element.
We have \cite{Zha98}
\bea
&&\L_{2N}=\frac{1}{2N}\sum_{k=1}^{2N}(-1)^{k+1}\ve_k,\no\\
&&\L_i=d+\sum_{k=1}^i(-1)^{k+1}\ve_k-\frac{i}{2N}\sum_{k=1}^{2N}
   (-1)^{k+1}\ve_k,
\eea
where $i=0,1,\cdots,2N-1$. 

The quantum affine superalgebra $\uqgnh$ is a quantum (or $q$-)
deformation of the universal enveloping algebra of $\hG$ and is
generated by the Chevalley generators $\{e_i,\;f_i\;q^{h_j},\;d|
i=0,1,\cdots,2N-1,\;j=0,1,\cdots,2N\}$.
The $\Z2$-grading of the Chevalley generators is $[e_i]=[f_i]=1,~
i=0,1,\cdots,2N-1$ and zero
otherwise. The defining relations are
\bea
&&hh'=h'h,~~~~~~\forall h\in \tilde{\cal H},\no\\
&&q^{h_j}e_iq^{-h_j}=q^{a_{ij}}e_i,~~~~[d, e_i]=\d_{i0}e_i,\no\\
&&q^{h_j}f_iq^{-h_j}=q^{-a_{ij}}f_i,~~~~[d,f_i]=-\d_{i0}f_i,\no\\
&&[e_i,f_{i'}]=\d_{ii'}\frac{q^{h_i}-q^{-h_i}}{q-q^{-1}},\no\\
&&[e_i,e_{i'}]=[f_i,f_{i'}]=0,~~~~{\rm for}~~a_{ii'}=0,\no\\
&&[[e_0,e_1]_{q^{-1}}, [e_0,e_{2N-1}]_q]=0,\no\\
&&[[e_l,e_{l-1}]_{q^{(-1)^l}}, [e_l,e_{l+1}]_{q^{(-1)^{l+1}}}]=0,\no\\
&&[[e_{2N-1},e_{2N-2}]_{q^{-1}}, [e_{2N-1},e_0]_q]=0,\no\\
&&[[f_0,f_1]_{q^{-1}}, [f_0,f_{2N-1}]_q]=0,\no\\
&&[[f_l,f_{l-1}]_{q^{(-1)^l}}, [f_l,f_{l+1}]_{q^{(-1)^{l+1}}}]=0,\no\\
&&[[f_{2N-1},f_{2N-2}]_{q^{-1}}, [f_{2N-1},f_0]_q]=0,~~~~l=1,2,\cdots,2N-2.
\eea
Here and throughout, $[a,b]_x\equiv ab-(-1)^{[a][b]}x ba$ 
and $[a,b]\equiv [a,b]_1$. The four-th order $q$-Serre relations are
obtained by using Yamane's Dynkin diagram procedure \cite{Yam96}.

$\uqgnh$ is a $\Z2$-graded quasi-triangular Hopf algebra endowed with
the following coproduct $\D$, counit $\e$ and antipode $S$:
\bea
\D(h)&=&h\otimes 1+1\otimes h,\no\\
\D(e_i)&=&e_i\otimes 1+q^{h_i}\otimes e_i,~~~~
\D(f_i)=f_i\otimes q^{-h_i}+1\otimes f_i,\no\\
\e(e_i)&=&\e(f_i)=\e(h)=0,\no\\
S(e_i)&=&-q^{-h_i}e_i,~~~~S(f_i)=-f_iq^{h_i},~~~~S(h)=-h,\label{e-s} 
\eea
where $i=0,1,\cdots,2N-1$ and $h\in \tilde{\cal H}$.
Notice that the antipode $S$ is a $\Z2$-graded algebra anti-homomorphism.
Namely, for any homogeneous elements $a,b\in\uqgnh$
$S(ab)=(-1)^{[a][b]}S(b)S(a)$, which extends to inhomogeneous elements
through linearity. Moreover,
\beq
S^2(a)=q^{-2\rho}\,a\,q^{2\rho},~~~~\forall a\in\uqgnh,
\eeq
where $\rho$ is an element in $\tilde{\H}$
such that $(\rho,\a_i)=(\a_i,\a_i)/2$ for any simple root
$\a_i,~i=0,1,2,\cdots,2N-1$. Explicitly,
\beq
\rho=\frac{1}{2}\sum_{k=1}^{2N}(-1)^k\ve_k,
\eeq
which coincides with $\bar{\rho}$, the half-sum of positive roots
of $\G$ in the present simple root system.
The multiplication rule on the tensor products is $\Z2$-graded:
$(a\otimes b)(a'\otimes b')=(-1)^{[b][a']}(aa'\otimes bb')$ for any
homogeneous elements $a,b,a',b'\in \uqgnh$.
We also introduce the element in $\tilde{\H}$
\beq
\tilde{\rho}=\sum_{i=0}^{2N-1}\L_i+\xi N\L_{2N},\label{rho-t}
\eeq
which gives the principal gradation
\beq
[\tilde{\rho},e_i]=e_i,~~~~[\tilde{\rho},f_i]=-f_i,~~~~i=0,1,\cdots,2N-1.
\eeq
In (\ref{rho-t}),  $\xi$ is an arbitrary constant.

$\uqgnh$ can also be realized in terms of the Drinfeld generators
\cite{Dri88} $\{X^{\pm,i}_m,\; H^j_n,\; q^{\pm H^j_0}$, $c,\; d |
m\in {\bf Z},\; n\in{\bf Z}-\{0\},\; i=1,2,\cdots,2N-1,\;
j=1,2,\cdots,2N\}$.  The $\Z2$-grading of the Drinfeld generators is
given by $[X^{\pm,i}_m]=1$ for all $i=1,\cdots,2N-1,\;m\in{\bf Z}$ and 
$[H^j_n]=[H^j_0]=[c]=[d]=0$ for all $j=1,\cdots,2N,\; n\in{\bf Z}-\{0\}$.
The relations 
satisfied by the Drinfeld generators read (see \cite{Yam96,Zha97}
for the Drinfeld realization of $\uqgnh$  in the distinguished
system of simple roots)
\bea
&&[c,a]=[d,H^j_0]=[H^j_0, H^{j'}_n]=0,~~~~\forall a\in\uqgnh\no\\
&&q^{H^j_0}X^{\pm, i}_nq^{-H^j_0}=q^{\pm a_{ij}}X^{\pm, i}_n,\no\\
&&[d,X^{\pm,i}_n]=nX^{\pm,i}_n,~~~[d,H^j_n]=nH^j_n,\no\\
&&[H^j_n, H^{j'}_m]=\d_{n+m,0}\frac{[a_{jj'}n]_q[nc]_q}{n},\no\\
&&[H^j_n,
   X^{\pm, i}_m]=\pm\frac{[a_{ij}n]_q}{n}X^{\pm,i}_{n+m}q^{\mp|n|c/2},\no\\
&&[X^{+,i}_n, X^{-,i'}_m]=\frac{\d_{ii'}}{q-q^{-1}}\lt(q^{\frac{c}{2}(n-m)}
  \psi^{+,i}_{n+m}-q^{-\frac{c}{2}(n-m)}\psi^{-,i}_{n+m}\rt),\no\\
&&[X^{\pm,i}_n, X^{\pm, i'}_m]=0,~~~~{\rm for}~a_{ii'}=0,\no\\
&&[X^{\pm,i}_{n+1}, X^{\pm,i'}_m]_{q^{\pm a_{ii'}}}
  -[X^{\pm,i'}_{m+1}, X^{\pm,i}_n]_{q^{\pm a_{ii'}}}=0,\no\\
&&[[X^{\pm,l}_m,X^{\pm,l-1}_{m'}]_{q^{(-1)^l}},[X^{\pm,l}_n,
  X^{\pm,l+1}_{n'}]_{q^{(-1)^{l+1}}}]\no\\
&&~~~+[[X^{\pm,l}_n,X^{\pm,l-1}_{m'}]_{q^{(-1)^l}},[X^{\pm,l}_m,
  X^{\pm,l+1}_{n'}]_{q^{(-1)^{l+1}}}]=0, ~~~l=2,\cdots,2N-2.
  \label{drinfeld}
\eea
where $[x]_q=(q^x-q^{-x})/(q-q^{-1})$ and
 $\psi^{\pm,j}_{n}$ are related to $H^j_{\pm n}$ by relations
\beq
\sum_{n\in{\bf Z}}\psi^{\pm,j}_{n}z^{- n}=q^{\pm H^j_0}\exp\lt(
  \pm(q-q^{-1})\sum_{n>0}H^j_{\pm n}z^{\mp n}\rt).
\eeq
The following relations can be proved by induction:
\bea
H^j_n&=&\frac{1}{q-q^{-1}}\sum_{p_1+2p_2+\cdots+np_n=n}\frac{(-1)^{\sum\,p_i-1}
       \;(\sum\,p_i-1)!}
       {p_1!\cdots p_n!}(q^{-H^j_0}\psi^{+,j}_1)^{p_1}\cdots 
       (q^{-H^j_0}\psi^{+,j}_n)^{p_n},\no\\
H^j_{-n}&=&\frac{1}{q^{-1}-q}\sum_{p_1+2p_2+\cdots+np_n=n}
      \frac{(-1)^{\sum\,p_i-1}\;(\sum\,p_i-1)!}
       {p_1!\cdots p_n!}(q^{H^j_0}\psi^{-,j}_{-1})^{p_1}\cdots 
       (q^{H^j_0}\psi^{-,j}_{-n})^{p_n}.
\eea

The Chevalley generators are related to the Drinfeld generators by the
formulae
\bea
&&h_i=H_0^i,~~~e_i=X^{+,i}_0,~~~f_i=X^{-,i}_0,~~~
  i=1,2,\cdots,2N-1,\no\\
&&h_{2N}=H^{2N}_0,~~~h_0=c-\sum_{k=1}^{2N-1}H^k_0,\no\\
&&e_0=[X^{-,2N-1}_0,[X^{-,2N-2}_0,\cdots,[X^{-,3}_0,[X^{-,2}_0,
  X^{-,1}_1]_{q}]_{q^{-1}}\cdots]_{q}]_{q^{-1}}
  \,q^{-\sum_{k=1}^{2N-1}H^k_0},\no\\
&&f_0=(-1)^Nq^{\sum_{k=1}^{2N-1}H^k_0}[[\cdots [[X^{+,1}_{-1},
  X^{+,2}_0]_{q^{-1}}
  X^{+,3}_0]_{q},\cdots,X^{+,2N-2}_0]_{q^{-1}},X^{+,2N-1}_0]_{q}.
\eea

The coproduct of the Drinfeld generators is not known in full. 
However, for our purpose it suffices to derive the ``main terms" \cite{Cha91}
in the coproduct formulae. We have
\begin{Proposition}\label{main-terms}: 
For $m\in {\bf Z}_{\geq 0}$, $n\in {\bf Z}_{>0}$
and $i=1,2,\cdots,2N-1$,
\bea
\D(X^{+,i}_m)&=&X^{+,i}_m\otimes q^{mc}+q^{H^i_0+2mc}\otimes
   X^{+,i}_m\no\\
& &+\sum_{k=0}^{m-1} q^{\frac{1}{2}(m+3k)c}\psi^{+,i}_{m-k}\otimes
   q^{(m-k)c}X^{+,i}_k~~~{\rm mod}~N_-\otimes N^2_+,\no\\
\D(X^{+,i}_{-n})&=&X^{+,i}_{-n}\otimes q^{-nc}+q^{-H^i_0}\otimes
   X^{+,i}_{-n}\no\\
& &+\sum_{k=1}^{n-1} q^{\frac{1}{2}(n-k)c}\psi^{-,i}_{k-n}\otimes
   q^{(k-n)c}X^{+,i}_{-n}~~~{\rm mod}~N_-\otimes N^2_+,\no\\
\D(X^{-,i}_{n})&=&X^{-,i}_{n}\otimes q^{H^i_0}+q^{nc}\otimes
   X^{-,i}_{n}\no\\
& &+\sum_{k=1}^{n-1} q^{(n-k)c}X^{-,i}_k\otimes 
   q^{\frac{1}{2}(k-n)c}\psi^{+,i}_{n-k}~~~{\rm mod}~N^2_-\otimes N_+,\no\\
\D(X^{-,i}_{-m})&=&X^{-,i}_{-m}\otimes q^{-H^i_0-2mc}+q^{-mc}\otimes
   X^{-,i}_{-m}\no\\
& &+\sum_{k=0}^{m-1} q^{(k-m)c}X^{-,i}_{-k}\otimes 
   q^{-\frac{1}{2}(m+3k)c}\psi^{-,i}_{k-m}~~~{\rm mod}~N^2_-\otimes N_+,\no\\
\D(H^i_n)&=&H^i_n\otimes q^{\frac{1}{2}nc}+q^{\frac{3}{2}nc}\otimes
    H^i_n~~~{\rm mod}~N_-\otimes N_+,\no\\
\D(H^i_{-n})&=&H^i_{-n}\otimes q^{-\frac{3}{2}nc}+q^{-\frac{1}{2}nc}\otimes
    H^i_{-n}~~~{\rm mod}~N_-\otimes N_+,
\eea
where $N_\pm$ and $N^2_\pm$ are the left ideals generated by $X^{\pm,k}_l$
and $X^{\pm,k}_lX^{\pm,k'}_{l'}$, $k,\;k'=1,\cdots,2N-1;\; 
l,\;l'\in {\bf Z}$, respectively.
\end{Proposition}

\noindent{\it Remark.} (i) We do not write down the formulae for
$\D(H^{2N}_{\pm n})$ because they are not needed in this paper. 
$\D(H^{2N}_{\pm n})$ can be determined by requiring that $\D$ preserves
the commutation relations (\ref{drinfeld}). (ii)
Modulo $N_+\otimes N_-+N_-\otimes N_+$, the elements
$\psi^{\pm,i}_{\pm n}~(n\geq 0)$ are group-like:
\bea
&&\D(\psi^{+,i}_n)=\sum_{k=0}^n q^{\frac{3}{2}kc}\psi^{+,i}_{n-k}
  \otimes q^{\frac{1}{2}(n-k)c}\psi^{+,i}_k,\no\\
&&\D(\psi^{-,i}_{-n})=\sum_{k=0}^n q^{-\frac{1}{2}kc}\psi^{-,i}_{k-n}
  \otimes q^{\frac{3}{2}(k-n)c}\psi^{-,i}_{-k}.
\eea

Define the Drinfeld currents or generating functions,
\beq
X^{\pm,i}(z)=\sum_{n\in{\bf Z}}X^{\pm,i}_n z^{-n-1},~~~~
\psi^{\pm,j}(z)=\sum_{n\in{\bf Z}}\psi^{\pm,j}_n z^{-n}
\eeq
In terms of these currents, (\ref{drinfeld}) read
\bea
&&\psi^{\pm,j}(z)\psi^{\pm,j'}(w)=\psi^{\pm,j'}(w)\psi^{\pm,j}(z),\no\\
&&\psi^{+,j}(z)\psi^{-,j'}(w)=\frac{(z-wq^{c+a_{jj'}})(z-wq^{-c-a_{jj'}})}
  {(z-wq^{c-a_{jj'}})(z-wq^{-c+a_{jj'}})}\psi^{-,j'}(w)
  \psi^{+,j}(z),\no\\
&&\psi^{+,j}(z)X^{\pm,i}(w)=q^{\pm a_{ij}}\frac{z-wq^{\mp\frac{c}{2}
  \mp a_{ij}}}{z-wq^{\mp\frac{c}{2}\pm a_{ij}}}
  X^{\pm,i}(w)\psi^{+,j}(z),\no\\
&&\psi^{-,j}(z)X^{\pm,i}(w)=q^{\pm a_{ij}}\frac{z-wq^{\pm\frac{c}{2}
  \mp a_{ij}}}{z-wq^{\pm\frac{c}{2}\pm a_{ij}}}
  X^{\pm,i}(w)\psi^{-,j}(z),\no\\
&&[X^{+,i}(z), X^{-,i'}(w)]=\frac{\d_{ii'}}{(q-q^{-1})zw}\lt(
  \d(\frac{w}{z}q^c)\psi^{+,i}(wq^{\frac{c}{2}})
  -\d(\frac{w}{z}q^{-c})\psi^{-,i}(wq^{-\frac{c}{2}})\rt),\no\\
&&X^{\pm,i}(z)X^{\pm,i'}(w)+X^{\pm,i'}(w)X^{\pm,i}(z)=0,~~~~
  {\rm for}~a_{ii'}=0,\no\\
&&(z-wq^{\pm a_{ii'}})X^{\pm,i}(z)X^{\pm,i'}(w)+
  (zq^{\pm a_{ii'}}-w)X^{\pm,i'}(w)X^{\pm,i}(z)=0,\no\\
&&\{[[X^{\pm,l}(z_1),X^{\pm,l-1}(z)]_{q^{(-1)^l}},[X^{\pm,l}(z_2),
  X^{\pm,l+1}(w)]_{q^{(-1)^{l+1}}}]\}\no\\
&&~~~~+\{z_1\leftrightarrow z_2\}=0,~~~~l=2,3,\cdots,2N-2.
  \label{current}
\eea
These current commutation relations can be derived from the super
version \cite{Zha97,Gou98} of the RS algebra \cite{Res90} by means of the Gauss
decomposition technique of Ding and Frenkel \cite{Din93}.

\sect{Level-Zero Representation}

We consider the evaluation representation $V_z$ of $\uqgnh$,
where $V$ is an $2N$-dimensional graded
vector space with basis vectors $\{v_1,v_2,\cdots,v_{2N}\}$. The
$\Z2$-grading of the basis vectors is chosen to be $[v_j]=\frac{(-1)^j
+1}{2}$. Let $e_{j,j'}$ be the $2N\times 2N$ matrix satisfying
$(e_{j,j'})_{kk'}=\d_{jk}\d_{j'k'}$ 
or equivalently $e_{i,j}v_k=\d_{jk}v_i$, (which
implies that for any operator $A$ its matrix elements $A_{j,i}$ are defined by
$Av_i=A_{j,i}v_j$). In the homogeneous gradation, the
Chevalley generators on 
$V_z$ are represented by 
\bea
&&e_i=e_{i,i+1},~~~f_i=(-1)^{i+1}e_{i+1,i},~~~i=1,2,\cdots,2N-1,\no\\
&&h_i=(-1)^{i+1}(e_{i,i}
  +e_{i+1,i+1}),~~~
  h_{2N}=\sum_{k=1}^{2N}(-1)^{k+1}e_{k,k},\no\\
&&e_0=ze_{2N,1},~~~f_0=-z^{-1}e_{1,2N},  ~~~h_0=-e_{1,1}-e_{2N,2N}.
\eea

Let $V^{*S}$ be the left dual module of $V$, defined by
\beq
(a\cdot v^*)(v)=(-1)^{[a][v^*]}v^*(S(a)v),~~~\forall a\in \uqgnh,\;
   v\in V,\; v^*\in V^*.
\eeq
Namely, the representations on
$V^{*S}$  are given by
\beq
\pi_{V^{*S}}(a)=\pi_V(S(a))^{st},~~~~\forall a\in\uqgnh,
\eeq
where $st$ denotes the supertansposition defined by
$(A_{i,j})^{st}=(-1)^{[j]([i]+[j])}A_{j,i}$. Note that in general
$((A_{i,j})^{st}))^{st}=(-1)^{[A]}A_{i,j}\neq A_{i,j}$. Let $V^{*S}_z$
be the $2N$-dimensional evaluation module corresponding to $V^{*S}$.
On $V^{*S}_z$, the Chevalley generators are represented by
\bea
&&e_i=-(-1)^i q^{(-1)^i}e_{i+1,i},~~~f_i=-q^{(-1)^{i+1}}e_{i,i+1},~~~
  i=1,2,\cdots,2N-1,\no\\
&&h_i=(-1)^i(e_{i,i}+e_{i+1,i+1}),~~~h_{2N}=\sum_{k=1}^{2N}
  (-1)^ke_{k,k},\no\\
&&e_0=zqe_{1,2N},~~~f_0=z^{-1}q^{-1}e_{2N,1},~~~
  h_0=e_{1,1}+e_{2N,2N}.
\eea

\begin{Proposition}\label{level-one}: The Drinfeld generators are represented
on $V_z$ by
\bea
H^i_m&=&(-1)^{i+1}\frac{[m]_q}{m}q^{(-1)^im}\lt(q^{x_i}z\rt)^m
  (e_{i,i}+e_{i+1,i+1}),\no\\
H^{2N}_m&=&z^m\frac{[2m]_q}{m}\lt[-q^m\sum_{l=1}^Ne_{2l,2l}\rt.\no\\
& &\lt. +\sum_{l=1}^N\lt(y+(l-1)(1-q^m)\rt)(e_{2l-1,2l-1}+
  e_{2l,2l})\rt],\no\\
H^i_0&=&(-1)^{i+1}(e_{i,i}+e_{i+1,i+1}),~~~~H^{2N}_0=\sum_{k=1}^{2N}
  (-1)^{k+1}e_{k,k},\no\\
X^{+,i}_m&=&\lt(q^{x_i}z\rt)^me_{i,i+1},~~~~
  X^{-,i}_m=(-1)^{i+1}\lt(q^{x_i}z\rt)^me_{i+1,i},
\eea
and on $V^{*S}_z$ by
\bea
H^i_m&=&(-1)^{i}\frac{[m]_q}{m}q^{(-1)^{i+1}m}\lt(q^{-x_i}z\rt)^m
  (e_{i,i}+e_{i+1,i+1}),\no\\
H^{2N}_m&=&-z^m\frac{[2m]_q}{m}\lt[-q^{-m}\sum_{l=1}^Ne_{2l,2l}\rt.\no\\
& &\lt. +\sum_{l=1}^N\lt(-y^*+(l-1)(1-q^{-m})\rt)(e_{2l-1,2l-1}+
  e_{2l,2l})\rt],\no\\
H^i_0&=&(-1)^{i}(e_{i,i}+e_{i+1,i+1}),~~~~H^{2N}_0=\sum_{k=1}^{2N}
  (-1)^{k}e_{k,k},\no\\
X^{+,i}_m&=&-(-1)^iq^{(-1)^i}\lt(q^{-x_i}z\rt)^me_{i+1,i},~~~~
  X^{-,i}_m=-q^{(-1)^{i+1}}\lt(q^{-x_i}z\rt)^me_{i,i+1},
\eea
where $i=1,\cdots,2N-1$,  $x_i=\sum_{l=1}^i (-1)^{l+1}=\frac{(-1)^{i+1}+1}{2}$
and $y,\;y^*$ are arbitrary constants.
\end{Proposition}

\sect{Free Boson Realization at Level One}

We use the notations similar to those in \cite{Awa94,Kim97}.
Let us introduce bosonic oscillators $\{A^j_n,\;c^l_n,\;Q_{A^j},\;
Q_{c^l}|n\in{\bf Z}, j=1,2,\cdots, 2N,\;l=1,2,\cdots,N\}$ which 
satisfy the commutation relations
\bea
&&[A^j_n, A^{j'}_m]=\d_{n+m,0}\frac{[a_{jj'}n]_q[n]_q}{n},~~~~~
  [A^j_0, Q_{A^{j'}}]=a_{jj'},\no\\
&&[c^l_n, c^{l'}_m]=\d_{ll'}\d_{n+m,0}\frac{[n]_q^2}{n},~~~~~
  [c^l_0, Q_{c^{l'}}]=\d_{ll'}.\label{oscilators}
\eea
The remaining commutation relations are zero. Introduce the currents
\bea
&&H^j(z;\k)=Q_{A^j}+A^j_0\ln z
   -\sum_{n\neq 0}\frac{A^j_n}{[n]_q}q^{\k |n|}z^{-n},\no\\
&&c^l(z)=Q_{c^l}+c^l_0\ln z-\sum_{n\neq 0}\frac{c^l_n}{[n]_q}
   z^{-n}
\eea
and  set
\bea
H^j_\pm(z)&=&H^j(q^{\pm\frac{1}{2}}z;-\frac{1}{2})-H^j(q^{\mp\frac{1}{2}}z;
  \frac{1}{2})\no\\
&=&\pm(q-q^{-1})\sum_{n>0}A^j_{\pm n}z^{\mp n}\pm A^i_0\ln q.
\eea
We make a basis transformation and express $A^j_n$ and $Q_{A^j}$ 
in terms of a new set of bosonic oscillators
$\{a^j_n,\;Q_{a^j}|j=1,2,\cdots,2N\}$ as
\bea
&&A^i_n=(-1)^{i+1}\lt(a^i_n+a^{i+1}_n\rt),~~~~
  A^{2N}_n=\frac{q^n+q^{-n}}{2}\sum_{l=1}^{2N} (-1)^{l+1}a^l_n,\no\\
&&Q_{A^i}=(-1)^{i+1}\lt(Q_{a^i}+Q_{a^{i+1}}\rt),~~~~
  Q_{A^{2N}}=\sum_{l=1}^{2N}(-1)^{l+1}Q_{a^l},\label{A-a}
\eea
where $i=1,\cdots,2N-1$ and
$\{a^j_n,\;Q_{a^j}\}$ satisfy the commutation relations
\beq
[a^j_n, a^{j'}_m]=(-1)^{j+1}\d_{jj'}\d_{n+m,0}\frac{[n]_q^2}{n},~~~~~
  [a^j_0, Q_{a^{j'}}]=(-1)^{j+1}\d_{jj'}.
\eeq
Now we state our main result in this section on the free boson
realization of $\uqgnh$ at level one.
\begin{Theorem}\label{free-boson}: 
The Drinfeld generators of $\uqgnh$  at level one are
realized by the free boson fields as 
\bea
&&c=1,\no\\
&&\psi^{\pm,j}(z)=e^{H^j_\pm(z)},~~~~j=1,2,\cdots,2N,\no\\
&&X^{\pm,i}(z)=:e^{\pm H^i(z;\mp\frac{1}{2})}\;Y^{\pm,i}(z):F^{\pm,i},~~~~
   i=1,2,\cdots,2N-1,
\eea
where 
\bea
&&F^{\pm,2k-1}=\prod_{l=1}^{k-1}e^{\pm\sqrt{-1}\pi a^{2l-1}_0},~~~~~
   F^{\pm,2k}=\prod_{l=1}^{k}e^{\mp\sqrt{-1}\pi a^{2l-1}_0},\no\\
&&Y^{+,2k-1}(z)=e^{c^k(z)},\no\\
&&Y^{-,2k-1}(z)=\frac{1}{z(q-q^{-1})}\lt(e^{-c^k(qz)}-
  e^{-c^k(q^{-1}z)}\rt),\no\\
&&Y^{+,2k}(z)=Y^{-,2k-1}(z)=\frac{1}{z(q-q^{-1})}\lt(e^{-c^k(qz)}-
  e^{-c^k(q^{-1}z)}\rt),\no\\
&&Y^{-,2k}(z)=-Y^{+,2k-1}(z)=-e^{c^k(z)},~~~~k=1,2,\cdots,N.
   \label{boson}
\eea
\end{Theorem}

\noindent{\it Proof.} We prove this theorem by checking that they
satisfy the defining relations (\ref{current}) of $\uqgnh$ with $c=1$.
It is easily seen that the first two relations in (\ref{current})
are true by construction.
The third and fourth ones follow from the definition
of $X^{\pm,i}(z)$ and the commutativity between $a^j_n$ and $c^l_n$. So
we only need to check the last three relations in (\ref{current}).

We write
\beq
Z^{\pm,i}(z)=:e^{\pm H^i(z;\mp\frac{1}{2})}: F^{\pm,i}.
\eeq
It is easily shown that
\bea
Z^{+,i}(z)Z^{+,i'}(w)&=&\lt\{
\begin{array}{l}
:Z^{+,i}(z)Z^{+,i'}(w):~~~~{\rm for}~~a_{ii'}=0~{\rm and}~ i\leq i',\\
-:Z^{+,i}(z)Z^{+,i'}(w):~~~~{\rm for}~~a_{ii'}=0~{\rm and}~ i>i',\\
(z-q^{-1}w)\;:Z^{+,i}(z)Z^{+,i'}(w):~~~~{\rm for}~~a_{ii'}=1~{\rm and}~ i<i',\\
-(z-q^{-1}w)\;:Z^{+,i}(z)Z^{+,i'}(w):~~~~{\rm for}~~a_{ii'}=1~{\rm
  and}~i>i',\\
(z-q^{-1}w)^{-1}\;:Z^{+,i}(z)Z^{+,i'}(w):~~~~{\rm for}~~a_{ii'}=-1,
\end{array}
\rt.\label{z+z+}\\
Z^{-,i}(z)Z^{-,i'}(w)&=&\lt\{
\begin{array}{l}
:Z^{-,i}(z)Z^{-,i'}(w):~~~~{\rm for}~~a_{ii'}=0~{\rm and}~ i\leq i',\\
-:Z^{-,i}(z)Z^{-,i'}(w):~~~~{\rm for}~~a_{ii'}=0~{\rm and}~ i>i'\\
(z-qw)\;:Z^{-,i}(z)Z^{-,i'}(w):~~~~{\rm for}~~a_{ii'}=1~{\rm and}~ i<i',\\
-(z-qw)\;:Z^{-,i}(z)Z^{-,i'}(w):~~~~{\rm for}~~a_{ii'}=1~{\rm and}~ i>i',\\
(z-qw)^{-1}\;:Z^{-,i}(z)Z^{-,i'}(w):~~~~{\rm for}~~a_{ii'}=-1,
\end{array}
\rt.\label{z-z-}\\
Z^{+,i}(z)Z^{-,i'}(w)&=&\lt\{
\begin{array}{l}
:Z^{+,i}(z)Z^{-,i'}(w):~~~~{\rm for}~~a_{ii'}=0~{\rm and}~ i\leq i',\\
-:Z^{+,i}(z)Z^{-,i'}(w):~~~~{\rm for}~~a_{ii'}=0~{\rm and}~ i>i',\\
(z-w)^{-1}\;:Z^{+,i}(z)Z^{-,i'}(w):~~~~{\rm for}~~a_{ii'}=1~{\rm and}~ i<i',\\
-(z-w)^{-1}\;:Z^{+,i}(z)Z^{-,i'}(w):~~~~{\rm for}~~a_{ii'}=1~{\rm and}~ i>i',\\
(z-w)\;:Z^{+,i}(z)Z^{-,i'}(w):~~~~{\rm for}~~a_{ii'}=-1,
\end{array}
\rt.\label{z+z-}
\eea
We have similar formulae for 
$Z^{+,i'}(w)Z^{+,i}(z),~
Z^{-,i'}(w)Z^{-,i}(z)$ and $Z^{-,i'}(w)Z^{+,i}(z)$. 

We now compute operator products $Y^{+,i}(z)Y^{+,i'}(w)$  and
$Y^{+,i}(z)Y^{-,i'}(w)$. It is easily seen from the definition of 
$Y^{\pm,i}(z)$ that the non-trivial products are those 
corresponding to $i=i'$ and $a_{ii'}=1$. 
Note that $a_{ii'}=1$ whenever $i=2k-1,~i'=2k$ (or $i=2k,~i'=2k-1$)
where $k=1,2,\cdots,N-1$. The corresponding operator products are
\bea
Y^{+,2k-1}(z)Y^{+,2k}(w)&=&\frac{1}{w(q-q^{-1})}\lt(
  \frac{:e^{c^k(z)}e^{-c^k(qw)}:}{z-qw}-\frac{:e^{c^k(z)}e^{-c^k(q^{-1}w)}:}
  {z-q^{-1}w}\rt),\no\\
Y^{-,2k-1}(z)Y^{-,2k}(w)&=&-\frac{1}{z(q-q^{-1})}\lt(
  \frac{:e^{-c^k(qz)}e^{c^k(w)}:}{qz-w}-\frac{:e^{-c^k(q^{-1}z)}e^{c^k(w)}:}
  {q^{-1}z-w}\rt),\no\\
Y^{+,2k-1}(z)Y^{-,2k}(w)&=&-(z-w)\;:Y^{+,2k-1}(z)Y^{-,2k}(w):\no\\
Y^{+,2k}(z)Y^{-,2k-1}(w)&=&\frac{1}{zw(q-q^{-1})}\lt(q(z-w)\;
  :e^{-c^k(qz)}e^{-c^k(qw)}:\rt.\no\\
& &  -q^{-1}(z-w)\;:e^{-c^k(q^{-1}z)}
  e^{-c^k(q^{-1}w)}:\no\\
& &-(qz-q^{-1}w)\;:e^{-c^k(qz)}e^{-c^k(q^{-1}w)}:\no\\
& &\lt.   -(q^{-1}z-qw)\;:e^{-c^k(q^{-1}z)}e^{-c^k(qw)}:\rt).\label{yy
   for a=1}
\eea
Since $Y^{+,2k}(z)=Y^{-,2k-1}(z)$ and $Y^{-,2k}(z)=-Y^{+,2k-1}(z)$, 
the products $Y^{+,i}(z)Y^{+,i}(w)$ and $Y^{\pm,i}(z)Y^{\mp,i}(w)$ can be
deduced from (\ref{yy for a=1}). For example,
\bea
Y^{+,2k}(z)Y^{-,2k}(w)&=&-Y^{+,2k}(z)Y^{+,2k-1}(w)\no\\
&=&\frac{1}{z(q-q^{-1})}\lt(
  \frac{:e^{-c^k(q^{-1}z)}e^{c^k(w)}:}{q^{-1}z-w}
  -\frac{:e^{-c^k(qz)}e^{c^k(w)}:}{qz-w}\rt).
\eea

By means of (\ref{z+z+}), (\ref{z+z-}), (\ref{z-z-}) and (\ref{yy for
a=1}) we can show that the last three relations in (\ref{current}) are
satisfied by (\ref{boson}). For instance,
\bea
[X^{+,2k}(z),X^{-,2k-1}(w)]&=&-\frac{1}{zw(q-q^{-1})^2}:Z^{+,2k}(z)
  Z^{-,2k-1}(w):\lt (\frac{1}{z-w}+\frac{1}{w-z}\rt)\no\\
& &\times\; \lt((qz-q^{-1}w)\;:e^{-c^k(qz)}e^{-c^k(q^{-1}w)}:\rt.\no\\
& &\lt.+(q^{-1}z-qw)\;:e^{-c^k(q^{-1}z)}e^{-c^k(qw)}:\rt)\no\\
&=&-\frac{1}{z^2w(q-q^{-1})^2}:Z^{+,2k}(z)
  Z^{-,2k-1}(w):\no\\
& &\times\; \d(\frac{w}{z})
   \lt((qz-q^{-1}w)\;:e^{-c^k(qz)}e^{-c^k(q^{-1}w)}:\rt.\no\\
& &\lt.+(q^{-1}z-qw)\;:e^{-c^k(q^{-1}z)}e^{-c^k(qw)}:\rt)=0.
\eea

\sect{Bosonization of Level-One Vertex Operators}

In this section, we study the level-one
vertex operators \cite{Fre92} of $\uqgnh$. Let
$V(\l)$ be the highest weight $\uqgnh$-module with the highest weight
$\l$. Consider the following intertwiners of $\uqgnh$-modules
\cite{Jim94}:
\bea
\Phi^{\mu V}_\l(z)&:&~~ V(\l)\longrightarrow V(\mu)\otimes V_z,
     \label{Phi}\\
\Phi^{\mu V^*}_\l(z)&:&~~ V(\l)\longrightarrow V(\mu)\otimes V^{*S}_z,
     \label{Phi*}\\
\Psi^{V\mu}_\l(z)&:&~~ V(\l)\longrightarrow V_z\otimes V(\mu),
     \label{Psi}\\
\Psi^{V^*\mu}_\l(z)&:&~~ V(\l)\longrightarrow V^{*S}_z\otimes V(\mu).
     \label{Psi*}
\eea
They are intertwiners in the sense that for any $x\in \uqgnh$
\beq
\Xi(z)\cdot x=\D(x)\cdot\Xi(z),~~~~\Xi(z)=\Phi^{\mu V}_\l(z),~
   \Phi^{\mu V^*}_\l(z),~\Psi^{V\mu}_\l(z),~\Psi^{V^*\mu}_\l(z).
   \label{intertwiner1}
\eeq
These intertwiners are even operators, that is their gradings are
$[\Phi^{\mu V}_\l(z)] = [\Phi^{\mu V^*}_\l(z)] = [\Psi^{V\mu}_\l(z)]
= [\Psi^{V^*\mu}_\l(z)] = 0$. According to \cite{Jim94},
$\Phi^{\mu V}_\l(z)~\lt(\Phi^{\mu V^*}_\l(z)\rt)$ is called type I
(dual) vertex operator and
$\Psi^{V\mu}_\l(z)~\lt(\Psi^{V^*\mu}_\l(z)\rt)$ type II (dual) vertex
operator.

We expand the vertex operators as a formal series \cite{Jim94}
\bea
&&\Phi^{\mu V}_\l(z)=\sum_{j=1}^{2N}\,\Phi^{\mu V}_{\l,j}(z)
   \otimes v_j, ~~~~
\Phi^{\mu V^*}_\l(z)=\sum_{j=1}^{2N}\,\Phi^{\mu V^*}_{\l,j}(z)
   \otimes v^*_j,\no\\ 
&&\Psi^{V\mu}_\l(z)=\sum_{j=1}^{2N}\,v_j\otimes \Psi^{V\mu}_{\l,j}(z),~~~~
\Psi^{V^*\mu}_\l(z)=\sum_{j=1}^{2N}\,v^*_j\otimes \Psi^{V^*\mu}_{\l,j}(z).
\eea
Then the intertwining property (\ref{intertwiner1}) reads in terms of
components
\bea
&&\sum\Phi^{\mu V}_{\l,j}(z)\,x\otimes v_j(-1)^{[v_j][x]}
  =\sum x_{(1)}\Phi^{\mu V}_{\l,j}(z)\otimes x_{(2)}v_j
  (-1)^{[v_j][x_{(2)}]},\no\\
&&\sum\Phi^{\mu V^*}_{\l,j}(z)\,x\otimes v^*_j(-1)^{[v^*_j][x]}
  =\sum x_{(1)}\Phi^{\mu V^*}_{\l,j}(z)\otimes x_{(2)}v^*_j
  (-1)^{[v^*_j][x_{(2)}]},\no\\
&&\sum v_j\otimes \Psi^{V\mu}_{\l,j}(z)x=\sum x_{(1)}v_j
   \otimes x_{(2)}\Psi^{V\mu}_{\l,j}(z)(-1)^{[v_j][x_{(2)}]},\no\\
&&\sum v^*_j\otimes \Psi^{V^*\mu}_{\l,j}(z)x=\sum x_{(1)}v^*_j
   \otimes x_{(2)}\Psi^{V^*\mu}_{\l,j}(z)(-1)^{[v^*_j][x_{(2)}]},
     \label{intertwiner2}
\eea
where we have used
the notation $\D(x)=\sum_x\,x_{(1)}\otimes x_{(2)}$ and
the fact that the vertex operators are even which implies $[\Phi^{\mu
V}_{\l,j}(z)] = [\Phi^{\mu V^*}_{\l,j}(z)] = [\Psi^{V\mu}_{\l,j}(z)]
= [\Psi^{V^*\mu}_{\l,j}(z)] = [v_j] = \frac{(-1)^j+1}{2}$.

Introduce the even  operators $\phi(z),\;\phi^*(z),\;\psi(z)$ and
$\psi^*(z)$,
\bea
&&\phi(z)=\sum_{j=1}^{2N}\,\phi_j(z)
   \otimes v_j, ~~~~
\phi^*(z)=\sum_{j=1}^{2N}\,\phi^*_j(z)
   \otimes v^*_j,\no\\ 
&&\psi(z)=\sum_{j=1}^{2N}\,v_j\otimes \psi_j(z),~~~~
\psi^*(z)=\sum_{j=1}^{2N}\,v^*_j\otimes \psi^*_j(z).
\eea
The grading of the components is given by $[\phi_j(z)] = [\phi^*_j(z)]
= [\psi_j(z)] = [\psi^*_j(z)] = \frac{(-1)^j+1}{2}$. Now we state
\begin{Proposition}\label{int-for-phipsi}: 
Assume that the operators  $\phi(z),\; \phi^*(z),\;
\psi(z),\; \psi^*(z)$ satisfy the 
intertwining relations (\ref{intertwiner2}). Then the operators $\phi(z)$
and $\psi(z)$ with respect to $V_z$ are determined by the components 
$\phi_{2N}(z)$ and $\psi_1(z)$, respectively.
With respect to $V^{*S}_z$, the operators $\phi^*(z)$ and $\psi^*(z)$ are
determined by
$\phi^*_1(z)$ and $\psi^*_{2N}(z)$, respectively. More explicitly, we have
for $l=1,2,\cdots, 2N-1$,
\bea
&&(-1)^l\,\phi_l(z)=[\phi_{l+1}(z), f_l]_{q^{(-1)^l}},\no\\
&&[\phi_l(z), f_l]_{q^{(-1)^l}}=0,\no\\
&&[\phi_k(z), f_l]=0,~~~~k\neq l,\,l+1.\label{phi-l-phi-2n}
\eea
\bea
&&q^{(-1)^{l+1}}\,\phi^*_{l+1}(z)=[\phi^*_{l}(z), f_l]_{q^{(-1)^{l+1}}},\no\\
&&[\phi^*_{l+1}(z), f_l]_{q^{(-1)^{l+1}}}=0,\no\\
&&[\phi^*_k(z), f_l]=0,~~~~k\neq l,\,l+1.\label{phi*-l-phi*-1}
\eea
\bea
&&\psi_{l+1}(z)=[\psi_{l}(z), e_l]_{q^{(-1)^{l+1}}},\no\\
&&[\psi_{l+1}(z), e_l]_{q^{(-1)^{l+1}}}=0,\no\\
&&[\psi_k(z), e_l]=0,~~~~k\neq l,\,l+1.\label{psi-l-psi-1}
\eea
\bea
&&(-1)^{l+1}\,q^{(-1)^{l}}\,\psi^*_{l}(z)=
  [\psi^*_{l+1}(z), e_l]_{q^{(-1)^{l}}},\no\\
&&[\psi^*_{l}(z), e_l]_{q^{(-1)^{l}}}=0,\no\\
&&[\psi^*_k(z), e_l]=0,~~~~k\neq l,\,l+1.\label{psi*-l-psi*-2n}
\eea
\end{Proposition}

Next we determine the relations of the components $\phi_{2N}(z),\;
\phi^*_1(z),\; \psi_1(z),\; \psi^*_{2N}(z)$ and the Drinfeld
generators. By means of proposition \ref{main-terms} and the
intertwining relations, we have
\begin{Proposition}:\label{rel-for-components} 
For $\phi(z)$ associated with $V_z$,
\bea
&&[\phi_{2N}(z), X^{+,i}(w)]=0,\no\\
&&q^{h_i}\phi_{2N}(z)q^{-h_i}=q^{-\d_{i,2N-1}}\phi_{2N}(z),\no\\
&&[H^i_n, \phi_{2N}(z)]=-\d_{i,2N-1}q^{\frac{3}{2}n}\frac{[n]_q}{n}
   z^n\phi_{2N}(z),\no\\
&&[H^i_{-n}, \phi_{2N}(z)]=-\d_{i,2N-1}q^{-\frac{1}{2}n}\frac{[n]_q}{n}
   z^{-n}\phi_{2N}(z);\label{phi-2n}
\eea
for $\phi^*(z)$ associated with $V^*_z$,
\bea
&&[\phi^*_1(z), X^{+,i}(w)]=0,\no\\
&&q^{h_i}\phi^*_1(z)q^{-h_i}=q^{\d_{i,1}}\phi^*_1(z),\no\\
&&[H^i_n, \phi^*_1(z)]=\d_{i,1}q^{\frac{3}{2}n}\frac{[n]_q}{n}
   z^n\phi^*_1(z),\no\\
&&[H^i_{-n}, \phi^*_1(z)]=\d_{i,1}q^{-\frac{1}{2}n}\frac{[n]_q}{n}
   z^{-n}\phi^*_1(z);\label{phi*-1}
\eea
for $\psi(z)$ associated with $V_z$,
\bea
&&[\psi_1(z), X^{-,i}(w)]=0,\no\\
&&q^{h_i}\psi_1(z)q^{-h_i}=q^{-\d_{i,1}}\psi_1(z),\no\\
&&[H^i_n, \psi_1(z)]=-\d_{i,1}q^{\frac{1}{2}n}\frac{[n]_q}{n}
   z^n\psi_1(z),\no\\
&&[H^i_{-n}, \psi_1(z)]=-\d_{i,1}q^{-\frac{3}{2}n}\frac{[n]_q}{n}
   z^{-n}\psi_1(z);\label{psi-1}
\eea
and for $\psi^*(z)$ associated with $V^*_z$,
\bea
&&[\psi^*_{2N}(z), X^{-,i}(w)]=0,\no\\
&&q^{h_i}\psi^*_{2N}(z)q^{-h_i}=q^{\d_{i,2N-1}}\psi^*_{2N}(z),\no\\
&&[H^i_n, \psi^*_{2N}(z)]=\d_{i,2N-1}q^{\frac{1}{2}n}\frac{[n]_q}{n}
   z^n\psi^*_{2N}(z),\no\\
&&[H^i_{-n}, \psi^*_{2N}(z)]=\d_{i,2N-1}q^{-\frac{3}{2}n}\frac{[n]_q}{n}
   z^{-n}\psi^*_{2N}(z).\label{psi*-2n}
\eea
\end{Proposition}

In order to obtain bosonized expressions of the vertex operators, we
introduce the following combinations of the Drinfeld generators:
\bea
&&A^{*i}_n=\sum_{l=1}^{2N-1}a^{-1}_{il}A^l_n+\frac{2}{q^n+q^{-n}}
  a^{-1}_{i,2N}A^{2N}_n,\no\\
&&A^{*i}_0=\sum_{l=1}^{2N}a^{-1}_{il}A^l_0,~~~~
  Q^*_{A^i}=\sum_{l=1}^{2N}a^{-1}_{il}Q_{A^l},~~~~i=1,2,\cdots,2N-1,
\eea
which satisfy the relations
\bea
&&[A^{*i}_n, A^{i'}_m]=\d_{ii'}\d_{n+m,0}\frac{[n]_q^2}{n},\no\\
&&[A^{*i}_n, A^{*i'}_m]=a^{-1}_{ii'}\,\d_{n+m,0}\frac{[n]^2_q}{n},\no\\
&&[A^{*i}_0, Q_{A^{i'}}]=\d_{ii'},~~~~
  [A^{i}_0, Q^*_{A^{i'}}]=\d_{ii'},\no\\
&&[A^{*i}_0, Q^*_{A^{i'}}]=a^{-1}_{ii'},~~~~
  i,i'=1,2,\cdots,2N-1.
\eea
Introduce the currents,
\beq
H^{*,j}(z;\k)=Q^*_{A^j}+A^{*j}_0\ln z-\sum_{n\neq 0}\frac{A^{*j}_n}{[n]_q}
   q^{k|n|}z^{-n}.
\eeq

Now we state our main theorem in this section on the bosonic realization of
the operators $\phi(z),\;\phi^*(z),\;\psi(z)$ and $\psi^*(z)$ at
level one. Thanks to the previous propositions,
we only need to determine one component for each operator and the other
components are represented by the integral of the currents.
\begin{Theorem}\label{boson-for-components}:  
The components $\phi_{2N}(z),\;\phi^*_1(z),\;
\psi_1(z)$ and $\psi^*_{2N}(z)$ can be realized explicitly as follows:
\bea
\phi_{2N}(z)&=&:e^{-H^{*,2N-1}(qz;\frac{1}{2})}e^{c^N(qz)}:
  e^{-\sqrt{-1}\pi a_0^1}
  \prod_{l=1}^{N-1} e^{-\sqrt{-1}\pi\frac{2N+l}{2N}a^{2l+1}_0},\no\\
\phi^*_1(z)&=&:e^{H^{*,1}(qz;\frac{1}{2})}:
  \prod_{l=1}^{N-1} e^{\sqrt{-1}\pi\frac{2N-l}{2N}
  a^{2l+1}_0},\no\\
\psi_1(z)&=&:e^{-H^{*,1}(qz;-\frac{1}{2})}:
  \prod_{l=1}^{N-1} e^{\sqrt{-1}\pi\frac{2N-l}{2N}
  a^{2l+1}_0},\no\\
\psi^*_{2N}(z)&=&\frac{1}{z(q-q^{-1})}:e^{H^{*,2N-1}(qz;-\frac{1}{2})}
  \lt(e^{-c^N(q^2z)}-e^{-c^N(z)}\rt):\no\\
& &\times\;  e^{-\sqrt{-1}\pi a_0^1}
  \prod_{l=1}^{N-1} e^{-\sqrt{-1}\pi\frac{2N+l}{2N}a^{2l+1}_0}.
\eea
\end{Theorem}

\noindent{\it Proof.} This theorem is proved by checking that the
construction satisfies all the intertwining relations.

\noindent{\it Remark.} The following inverse elements of the
extended Cartan matrix are needed to determine the cocycle factors
appearing in above theorem:
\bea
&&a^{-1}_{2N-1,2l}=a^{-1}_{2N-1,2l+1}=-\frac{l}{N},~~~~l=1,2,\cdots,N-1,\no\\
&&a^{-1}_{2N-1,1}=0,~~~~a^{-1}_{2N-1,2N}=\frac{1}{2N},\no\\
&&a^{-1}_{1,2l-1}=a^{-1}_{1,2l}=\frac{N-l}{N},~~~~l=1,2,\cdots,N-1,\no\\
&&a^{-1}_{1,2N-1}=0,~~~~a^{-1}_{1,2N}=\frac{1}{2N}.
\eea

We are now in a position to state the following result:
\begin{Proposition}:
The vertex operators $\Phi^{\mu V}_\l(z),\;
\Phi^{\mu V^*}_\l(z),\; \Psi^{V\mu}_\l(z)$ and $\Psi^{V^*\mu}_\l(z)$, if
they exist,
satisfy the same relations as the operators $\phi(z),\;\phi^*(z),\; \psi(z)$
and $\psi^*(z)$, respectively. 
\end{Proposition}
This proposition follows immediately from the fact that the formers
and the latters obey the same intertwining properties. 
Identifying $\Phi^{\mu V}_\l(z),\;
\Phi^{\mu V^*}_\l(z),\; \Psi^{V\mu}_\l(z)$ and $\Psi^{V^*\mu}_\l(z)$
with $\phi(z),\;\phi^*(z),\; \psi(z)$
and $\psi^*(z)$, respectively, then
the bosonic realization of the vertex operators is easily seen to be
given by propositions \ref{int-for-phipsi}, \ref{rel-for-components}
and theorem \ref{boson-for-components}.

\vskip.3in
\noindent {\bf Acknowledgements.}
The financial support from Australian Research 
Council through a Queen Elizabeth II Fellowship Grant is
gratefully acknowledged.


\end{document}